\documentclass[12pt,oneside]{article}
\usepackage{amsmath,amsthm,amsfonts,amssymb,MnSymbol,mathrsfs,latexsym}
\usepackage{mathtools}
\usepackage[english]{babel}

\usepackage{blindtext}
\usepackage{geometry}

\newcommand{\const}{\rm const}

  \newcommand{\Dom}{\rm  Dom}

\theoremstyle{plain}
\newtheorem{theorem}{Theorem}[section]

\newtheorem{proposition}[theorem]{Proposition}
\newtheorem{definition}{Definition}[section]
\newtheorem{remark}{Remark}[section]

\title{\large \textbf{Generalized Grand Lebesgue Spaces norm\\ estimations for some operators}}

\footnotesize\date{}

\author{\normalsize Maria Rosaria Formica ${}^1$,   \normalsize Eugeny Ostrovsky
${}^2$ and \normalsize Leonid Sirota ${}^3$}

\begin{document}

\maketitle

\begin{center}
{\footnotesize ${}^{1}$ Universit\`{a} degli Studi di Napoli \lq\lq Parthenope\rq\rq, via Generale Parisi 13,\\
Palazzo Pacanowsky, 80132,
Napoli, Italy.} \\

\vspace{2mm}

{\footnotesize e-mail: mara.formica@uniparthenope.it} \\

\vspace{4mm}

{\footnotesize ${}^{2,\, 3}$  Bar-Ilan University, Department of Mathematics and Statistics, \\
52900, Ramat Gan, Israel.} \\

\vspace{2mm}

{\footnotesize e-mail: eugostrovsky@list.ru}\\

\vspace{2mm}

{\footnotesize e-mail: sirota3@bezeqint.net} \\

\end{center}

\begin{abstract}
 We study moment rearrangement invariant spaces, which contain as particular cases the generalized Grand Lebesgue Spaces, and provide norm estimates for some operators, not necessarily linear, acting between some measurable rearrangement invariant spaces. \par
  \ The estimations are formulated in the terms of fundamental functions for these spaces.
\end{abstract}

\vspace{5mm}


\vspace{4mm}

\section{Notations. Statement of problem. Assumptions.}

\vspace{4mm}

 \hspace{3mm} Let $ \ (X = \{x\}, \cal{B},\mu)   \ $ and $ \ (Y = \{y\}, \cal{F}, \nu) \ $ be two measurable spaces equipped with
 non-trivial sigma-finite measures $ \ \mu, \ \nu, \ $ respectively
 and let $ \ T = \{t\} \subset \mathbb R_+ $ be a numerical (parametrical) set.
 \ Let $ \ Q \ $ be an operator, non necessarily linear, which maps numerical values
 functions $f$ defined on the set $ X \otimes T$, measurable with respect the variable $x\in X$ for all the values $ \ t \in T \
 $, into ones defined on the space $Y \otimes T$:
 $$
  u = Q[f], \ \ \ f: \ X \otimes T \to \mathbb R, \ \  u: \ Y \otimes T \to \mathbb R.
 $$
  \ It will be presumed that the function $ \ u = u(y,t) \ $,  for all the values $  t \in T \ $, is also measurable
  with respect to the variable $ y \in Y. \ $  \par

 \ For instance, an integral operator of the form

 $$
 u_0(y,t) := Q_0[f](y,t) = \int_X K(t, y,x, f(x,t)) \ \mu(dx),
 $$
of course, under appropriate natural conditions on the kernel $ \ K(\cdot). \ $ \par

\vspace{3mm}

 \ Denote as usually the classical Lebesgue - Riesz norms, more exactly a {\it family } of norms
 parametrized by means of the set $ \ T \ $
 of the functions $ \ f,u \ $ relative to the measurable spaces  $ \ (X, \mu) \ $ and $  \ (Y,\nu)$, respectively, as follows

$$
||f||_q = ||f||_q(t) \stackrel{def}{=} \left[ \ \int_X |f(x,t)|^q \
\mu(dx) \ \right]^{1/q}, \ q \in [1,\infty)
$$

$$
||u||_p = ||u||_p(t) \stackrel{def}{=} \left[ \ \int_Y |u(y,t)|^p \
\nu(dy) \ \right]^{1/p}, \ p \in [1,\infty).
$$

 \vspace{4mm}

 {\bf Our assumptions:} there exist pairs of numbers  $ \ (a,b), \ (c,d), \ $  such that
 $$
 1 \le a < b \le \infty, \hspace{3mm} 1 \le c < d \le \infty
 $$
 and there is a non-negative constant
 $ \ C = C(a,b; c,d)$ \  such that, for the set $ \ T = \{t\} \subset \mathbb R_+$,

\begin{equation} \label{assumption}
||u||_p \le C(a,b,c,d) \  t^{1/p - 1/q} \ ||f||_q, \hspace{3mm}
\forall t \in T,  \hspace{3mm} q \in (a,b),\  \ p \in (c,d).
\end{equation}

\vspace{4mm}

  \ Let us mention some examples of such  operators. These inequalities appear, e.g., in the article \cite{Mascia Nguyen}
devoted to the theory of Partial Differential Equations (PDEs) of
hyperbolic type. The applications of the theory of relations
\eqref{assumption} in the approximation theory, more precisely, in
the so-called Nikolskii inequality, may be found in \cite{Ostr
Nikol}.

\vspace{4mm}

\hspace{2mm} {\bf Our aim in this short report is to extend the last
estimate} \eqref{assumption} {\bf  and similar ones in the setting
of the Grand Lebesgue Spaces (GLS) instead of the classical Lebesgue
- Riesz spaces.} \par

\vspace{3mm}

\begin{remark}
{\rm  The case considered in \eqref{assumption} as well as its
generalizations cannot be deduced, in general case, from the one
considered in \cite{Ostrovsky5}.
 }
 \end{remark}

\vspace{4mm}

\begin{center}

 \ {\it \ Grand Lebesgue Spaces (GLS).}  \par

\end{center}

\vspace{4mm}

\ We recall here briefly some known definitions and facts from the
theory of Grand Lebesgue Spaces (GLS).
    \ Let $ 1 \le  a < b \le \infty$ and $ \ \psi = \psi(p), \ p \in (a,b)$, be a {\it strictly positive} measurable numerical valued
function, not necessary finite in $a,b$, such that $ \displaystyle
\inf_{p \in (a,b)} \psi(p) > 0$. \\
We write $ (a,b) := \Dom [\psi]$.

We denote with $ G \Psi(a,b)$ the set of all such functions
$\psi(p), \ p \in (a,b)$, for some  $ 1 \le a < b \le \infty \ $ and
put
$$
G\Psi := \displaystyle\bigcup_{1 \le a < b \le \infty} G \Psi(a,b).
$$
For instance, if $m = {\const} > 0$,
$$
\psi_m(p) := p^{1/m},  \ \ p \in [1,\infty)
$$
or, for $ 1 \le a < b < \infty, \ \alpha,\beta = \const \ge 0$,
$$
   \psi^{(a,b; \alpha,\beta)}(p) := (p-a)^{-\alpha} \ (b-p)^{-\beta}, \  \ p \in
   (a,b),
$$
the above functions belong to $ G \Psi(a,b)$.

The (Banach) Grand Lebesgue Space$ \ G \psi  = G\psi(a,b)$ consists
of all the real (or complex) numerical valued measurable functions
$f: \Omega \to \mathbb R$ having finite norm
\begin{equation} \label{norm psi}
    || \ f \ || = ||f||_{G\psi} \stackrel{def}{=} \sup_{p \in (a,b)} \left[ \frac{||f||_p}{\psi(p)} \right].
 \end{equation}
We write $ \ G\psi \ $ when $ \ a = 1\ $ and $ \ b= \infty. \ $ The
function $ \  \psi = \psi(p) \  $ is named the {\it generating
function } for the space $G \psi$ . \par If for instance
\begin{equation} \label{extremal case}
  \psi(p) = \psi^{(r)}(p) = 1, \ \  p = r;  \  \ \psi^{(r)}(p) = +\infty,   \ \ p \ne r,
\end{equation}
 where $ \ r = {\const} \in [1,\infty),  \ C/\infty := 0, \ C \in \mathbb R, \ $ (an extremal case), then the correspondent
 $ \  G\psi^{(r)}(p)  \  $ space coincides  with the classical Lebesgue - Riesz space $ \ L_r = L_r(\Omega). \ $ \par

\vspace{4mm}

These spaces, their particular cases and also generalizations of
them, have been widely investigated (see, e.g.,
\cite{Fiorenza2000,KozOs,Lifl,Ostrovsky2,Ostrovsky4} and references
therein). They play an important role in the theory of Partial
Differential Equations (PDEs) (see, e.g.,
\cite{AFFGR2020,Fiorenza-Formica-Gogatishvili-DEA2018,Fiorenza-Formica-Rakotoson-DIE2017,Greco-Iwaniec-Sbordone-1997}),
in interpolation theory (see, e.g.,
\cite{Fiorenza-Karadzhov-2004,fioforgogakoparakoNA}), in the theory
of Probability (see, e.g., \cite{Ermakov etc.
1986,KozOsSir2017,Ostrovsky3,Ostrovsky5,ForKozOstr_Lithuanian}), in
Statistics and in the theory of random fields (see, e.g., \cite{Buld
Koz AMS,KozOs,Ostrovsky4}, \cite[chapter 5]{Ostrovsky1}), in
Functional Analysis and so on.
\par

 These spaces are rearrangement invariant (r.i.) Banach function spaces; their fundamental function has been studied in  \cite{Ostrovsky4}. They
 do not coincide, in the general case, with the classical rearrangement invariant spaces: Orlicz, Lorentz, Marcinkiewicz, etc.
 (see \cite{Lifl,Ostrovsky2}).

\vspace{4mm}

 The belonging of a function $ f: \Omega \to \mathbb{R}$
to some $ G\psi$ space is closely related to its tail function
behavior
$$
 T_f(t) \stackrel{def}{=} {\bf P}(|f| \ge t), \ \ t \ge 0,
 $$
as $ \ t \to 0+ \ $ as well as when $ \ t \to \infty $ (see
\cite{KozOs,KozOsSir2017,KozOsSir2019}).

 \vspace{3mm}

 \hspace{3mm}  The so-called {\it fundamental function} $ \ \phi[G \psi](\delta), \ \delta \ge 0$, for these spaces has the form

\begin{equation} \label{fun funct}
\phi[G \psi](\delta) \stackrel{def}{=} \sup_{p \in Dom(\psi)} \left[ \ \frac{\delta^{1/p}}{\psi(p)} \ \right].
\end{equation}

 \vspace{3mm}

 \ In particular, if $ \ X = \mathbb R^d \ $ equipped with ordinary Lebesgue measure $ \ d\mu = dx, \ $ then

$$
\phi[L_r(\mathbb R^d)](\delta) = \delta^{1/r}, \ \delta \in [0,
\infty.]
$$

 \vspace{3mm}

  \hspace{3mm} This function plays a very important role in functional analysis, in particular in the theory of Fourier series
 \cite{Bennet Sharpley}, theory of probability, especially in the theory of random fields \cite{Ostrovsky4} etc.
 The fundamental function  $ \ \phi[G \psi](\delta)$ of the Grand Lebesgue Space is studied in details in
\cite{Ostrovsky4}, where many examples of these functions are
provided for different Grand Lebesgue Spaces. \par

\vspace{4mm}

\section{Main result.}

\vspace{4mm}

 \hspace{3mm} Let us assume that inequality \eqref{assumption} holds. We suppose now that the function $ \ f =
 f(x,t)$, for all the values $ \ t \in T \ $ and with respect to the variable $ \ x \in X $, belongs to some Grand Lebesgue Space
 $ \ G\psi_{a,b}: \ $
$$
||f||_{G\psi_{a,b}} < \infty, \ \ \ \ \forall t \in T
$$
and analogously assume that
$$
 ||u||_{G\nu_{c,d}} < \infty, \ \ \ \ \forall t \in T.
$$

 \ One can suppose, without loss of generality, $ \ ||u||_{G\nu} = ||f||_{G\psi} = 1. \ $ It follows from the
 direct definition of the norm in the GLS

\begin{equation} \label{key relat}
||u||_p \le \nu(p), \ p \in (c,d); \hspace{4mm} ||f||_q \le \psi(q),
\ q \in (a,b).
\end{equation}

\ We conclude from the assumption \eqref{assumption}

\begin{equation} \label{frac ineq}
||u||_p \cdot \left[ \ \frac{t^{1/q}}{\psi(q)} \ \right] \le C(a,b,c,d) \ ||f||_q \cdot \left[ \ \frac{t^{1/p}}{\nu(p)}
 \ \right], \  \ t \in T.
\end{equation}
 \ One can take the maximum over $ \ q \ $ and next with respect to the variable $ \ p \ $; using the direct definition of the
 fundamental function we deduce the following result.\par

\vspace{3mm}

\begin{proposition}\label{prop main}
\begin{equation} \label{main res}
\phi[G\psi](t) \cdot ||u||_{G\nu} \le C(a,b,c,d) \ \phi[G\nu](t)
\cdot ||f||_{G\psi},
\end{equation}
or equivalently
\begin{equation} \label{main refin}
\frac{||u||_{G\nu}}{\phi[G\nu](t)} \le C(a,b,c,d) \cdot
\frac{||f||_{G\psi}}{\phi[G\psi](t)}, \ \  t \in T.
\end{equation}

\end{proposition}

\vspace{4mm}

\begin{remark}
{\rm  It is easy to see that the last inequality is the direct
generalization of \eqref{assumption}; see \eqref{extremal case}. }
\end{remark}

\begin{remark}
{\rm  The last estimate is essentially non - improvable in the
general case, see e.g. \cite{Ostr Nikol}.}
\end{remark}

\vspace{4mm}

 \section{Some generalization.}

\vspace{5mm}

\begin{center}

{\it  Moment rearrangement invariant spaces.}

\end{center}

 \hspace{3mm} Let $ \  (Z, < \ \cdot \ > Z)  \ $ be another rearrangement invariant space, where  $ \ Z \ $  is a linear subset
 of the space of all measurable functions, mapping  $ \ X \to \mathbb R, \ $  (or alike $ \ Y \to \mathbb R. ) \ $ \par

\vspace{3mm}

\begin{definition} {\rm (see \cite{Ostr Nikol})}

{\rm  We will say that the space $ \ Z \ $ with the norm $ \ < \cdot
> Z \ $ is a {\it moment rearrangement invariant space,} briefly: {\it
m.r.i.}  space, and write $ \ Z = (Z,  < \cdot > Z) \  \in \ m.r.i.,
\ $ iff there exist $ \ \alpha,  \beta$, where \ $1 \le \alpha <
\beta \le \infty, \ $ and some {\it rearrangement invariant norm} $
\ < \cdot >Z \ $ defined on the space of real valued functions
defined on the interval $ \ (\alpha, \beta), \ $ non necessary to be
finite on all the functions, such that
\begin{equation} \label{mri space}
\forall f: X \to R \ \Rightarrow \ |||f|||_Z \stackrel{def}{=}
<h(\cdot)>Z,
\end{equation}
where in turn
$$
 h(p) = ||f||_p = ||f||_{p,X}.
$$
}
\end{definition}

 \vspace{3mm}

For instance, the mentioned before GLS  are m.r.i. spaces.

\noindent There are many examples of these spaces in \cite{Ostr
Nikol}, as well as some its generalizations and counterexamples. See
 also \cite{Astashkin}. \par

\vspace{3mm}

 \ Let us  describe the fundamental function for these spaces. Let $ \ \delta = \const \ge 0; \ $
 assume that there exists a measurable set $ \ A \ $ in the space $ \ Z \ $ for which $ \ \mu(A) = \delta. \ $
 Denote the set of all such values $ \ \delta \ $ as $ \Delta[Z] := \{\delta\}. \ $ \par
  The function $ \ \delta \to \ \kappa[Z](\delta)$, 
  represents the {\it fundamental function} for the space $ \ Z$, where

\begin{equation} \label{kappa fun}
\kappa[Z](\delta) := < g_{\delta}(\cdot)> Z,
\end{equation}
and
\begin{equation} \label{gdelta}
 g_{\delta}(p) \stackrel{def}{=}  \delta^{1/p}, \ \ \delta \in \Delta[Z].
\end{equation}

\vspace{3mm}

 \ In the role of the space $ \ Z \ $ will act both the spaces $ \ (X,Y). \ $  From \eqref{assumption} we have, as
 before,

\vspace{3mm}

\begin{equation} \label{mri assumption}
\forall t \in T \ \Rightarrow \  ||u||_p \cdot t^{1/q} \le
  C(a,b,c,d) \  t^{1/p} \ ||f||_q, \hspace{3mm} q \in (a,b); \ p \in (c,d).
\end{equation}

\vspace{3mm}

 \ We take in (\ref{mri assumption}) consequently the norm $ \ < \cdot >X \ $  and  $ \ < \cdot >Y \ $ from both the sides

\begin{equation} \label{prelim}
< \ u \ > Y \cdot \kappa[X](t) \le   C(a,b,c,d) \cdot < \ f \ > X \cdot \kappa[Y](t).
\end{equation}

\vspace{4mm}

To summarize, we conclude:

\begin{proposition}\label{prop main generalized}

Under the formulated notations and restrictions, we have

\begin{equation} \label{main mri}
\frac{< u > Y }{\kappa[Y](t)} \le C(a,b,c,d) \cdot \frac{< f > X }{\kappa[X](t)}, \hspace{3mm} t \in \Delta[X] \cap \Delta[Y],
\end{equation}
which is a direct generalization of the statement of Proposition
\ref{prop main}.
\end{proposition}

\vspace{4mm}

\begin{center}

{\it   More general source datum.}

\end{center}

\vspace{4mm}

 \hspace{3mm} Let us consider the following generalization of the relation \eqref{assumption}: the set $ \ T  = \{t\} \ $ may be
 arbitrary; suppose now that there exists two  {\it positive}  numerical valued   finite functions  $ \  A = A(t), \ B = B(t)  \ $
such that there is a finite nonnegative constant  $ \ D = D(a,b,c,d)
\ $ such that, for the set  $ \ T = \{t\} $,

\vspace{3mm}

\begin{equation} \label{general assumption}
\forall t \in T \ \Rightarrow \  ||u||_p \le D(a,b,c,d) \  A(t)^{1/p} \ B(t)^{- 1/q} \ ||f||_q, \hspace{3mm} q \in (a,b); \ p \in (c,d).
\end{equation}

\vspace{4mm}

We conclude under formulated notations and restrictions as in
Proposition \ref{prop main generalized}:

\begin{proposition}
\begin{equation} \label{main mri}
\frac{< u > Y }{\kappa[B(t)]} \le D(a,b,c,d)  \cdot \frac{< f > X }{\kappa[A(t)]}.
\end{equation}
\end{proposition}

\vspace{5mm}

 \section{Concluding remarks.}

\vspace{5mm}

{\bf A.} Open question: for which operators and pairs of r.i. spaces
holds true the last result \eqref{main mri}?

\vspace{3mm}

\noindent {\bf B.} It is interesting to find the exact exponential
{\it lower} norms estimations in the considered case.

\vspace{3mm}

\noindent {\bf C.} You also need to find some other spaces and
operators having the described features.

\vspace{6mm}

\emph{\textbf{Acknowledgements}.} {\footnotesize 
M.R. Formica is member of Gruppo Nazionale per l'Analisi Matematica,
la Probabilit\`{a} e le loro Applicazioni (GNAMPA) of the Istituto
Nazionale di Alta Matematica (INdAM) and member of the UMI group
\lq\lq Teoria dell'Approssimazione e Applicazioni (T.A.A.)\rq\rq }

\end{document}